\documentclass[3p]{elsarticle}
\usepackage{epsfig}
\newtheorem{remark}{Remark}[section]

\newtheorem{theorem}{Theorem}[section]

\newtheorem{example}{Example}[section]

\newtheorem{lemma}{Lemma}[section]
\newtheorem{Ack}{Acknowledgement}[section]

\newproof{pf}{Proof}
\newproof{pot}{Proof of Theorem}

\usepackage{amssymb}

\begin{document}

\begin{frontmatter}

\title{\bf Analysis of  biological integrate-and-fire oscillators}

\author{Marat Akhmet} \ead{marat@metu.edu.tr}

\address{Department of Mathematics, Middle East Technical University, 06531 Ankara, Turkey}

\begin{abstract}
We consider discontinuous dynamics of  integrate-and-fire models,  which  consist of pulse-coupled biological oscillators. A  thoroughly  constructed   map is in the basis of the analysis. Synchronization of non-identical oscillators is investigated.   Significant advances for the solution  of  second Peskin's conjecture have been made. Examples with numerical simulations are given to validate the theoretical results. Perspectives  are discussed.
\end{abstract}

\end{frontmatter}
\section{Introduction} 
The collective behavior of biological and chemical oscillators is a fascinating  topic that has  attracted a lot  of attention in the last  $50$ years \cite{buck1}-\cite{winfree1}.  An exceptional place in the  analysis belongs to  synchronization,  which  in its general sense is understood as phase locking, frequency  locking, and synchrony  itself, that is motion in unison \cite{buck1,ermen},\cite{ermen2}-\cite{strogatz2}. The  integrate-and-fire model of  the cardiac pacemaker  \cite{knight} was  developed by   C. Peskin \cite{p} to   a  population   of identical pulse-coupled  oscillators. It was  conjectured that  the model   self-synchronizes such  that: 
 \begin{itemize}
\item[(C1)]   For arbitrary initial conditions, the system approaches a state in which all the oscillators are firing synchronously. 
\item[(C2)]   This remains true even when the oscillators are not quite identical.
\end{itemize} 
The  conjecture  $(C1)$   is  solved in  \cite{p} for a  system with two oscillators, and in \cite{ms} for the generalized model of  two and more  oscillators. The last  paper  gave  start  to an intensive and productive investigation of the problem  and  its  applications \cite{epg,gerst,strogatz}, \cite{bot}-\cite{timme2}.  As far as we know  the conjecture $(C2)$  remains  unsolved. Even  a developed  non-identity concept  has  not been found in the literature.  

In the present paper we generalize the model, and propose  a version of  non-identity.  The model is considered such  that  perturbations save the synchronization.  These oscillators  are not  only  pulse-coupled, but  connected during the time between  moments of firing. That  is, the  modeling differential equations are  not  separated.  One can  see that  this approach  may   provide more biological  sense  to this theory.  The paper consists of main results, simulations and  discussion of  possible extensions. We believe that the results and proposals of this  paper  reveal  new perspectives  on the  study  of integrate-and-fire  models  of oscillators.  The main role  in our  analysis  is played  by  a  specially defined map. It is not  a Poincar\'e map, since  it transforms  the   coordinate of one oscillator to that  of another, and  the two  interchange  roles in the course of  the mapping. If there are more than two oscillators,  they  are used  in pairs to shape the  map with the interference of other oscillators  acting as  perturbation.  
   
  The main object of our investigation is an integrate-and-fire model, which consists of  $n$ non-identical  pulse-coupled  oscillators, $x_i, i =1,2,\ldots,n.$    If  the system does not fire the oscillators  satisfy  the following equations 
 \begin{eqnarray}
 \label{1}
 && 	x'_i =  f(x_i)  + \phi_i(x).
 \end{eqnarray}
 The domain  consists of  all points $x  = (x_1,x_2,\ldots,x_n)$ such  that $0 \le x_i  \le 1 +\zeta_i(x)$ for all $i =1,2,\ldots,n.$ When  the oscillator $x_j$ increases  from zero, and  meets the surface   such  that   $x_j(t) = 1 + \zeta_j(x(t)),$ then it  fires,   $x_j(t+) = 0.$ This firing changes the values of all oscillators with $i \not =j,$  
	 \begin{eqnarray}\label{2}
x_i(t+)=
\left\{ 
\begin{array}{ll} 
0,\:\: {\rm if} \:\:x_i(t) + \epsilon +\epsilon_{i} \ge 1+\zeta_i(x), \\ 
x_i(t) + \epsilon +\epsilon_{i},\:\:{\rm otherwise}. \end{array}
\right.
\end{eqnarray}
Thus, it  is assumed that if   $x_i(t)  \ge 1+\zeta_i(x) - \epsilon -\epsilon_{i},$ then the oscillator  fires, too.
    It is assumed also  that  there exist positive  constants $\mu_i$ and $\xi_i$  such  that  $|\phi_i(x)| < \mu_i$  and $|\zeta_i(x)| < \xi_i,$ for all $x$ and $i.$ In what  follows, we call the real numbers  $\epsilon,\mu_i, \xi_i, \epsilon_{i},$ {\it parameters}, assuming  the  first one is positive.  Moreover,  constants $\xi_i,\epsilon_{i},\mu_i,$  will be called {\it parameters of perturbation}. If all  of them  are   zeros, then  the model of identical oscillators is  obtained.  We assume that  $\epsilon + \epsilon_{i}> 0.$  That  is, an exhibitory  model  is under discussion.    The function  $f$  is positive valued and  lipschitzian.   Moreover, assume that  $\zeta_i$   are continuous and $\phi_i$ are locally  lipschitzian  for all $i.$
  
   The coupling  in the model is all-to-all such  that each   firing elicits  jumps in   all non-firing oscillators.   If several oscillators fire simultaneously, then  other oscillators react  as if just  one oscillator fires. In other words, any  firing  acts  only  as a signal  which abruptly  provokes  a   change of state. The intensity  of the signal is not important,  and pulse strengths are not additive.  A system of oscillators is synchronized if all of them fire in unison. 
    
     In the present analysis  we address   synchronization as well as the existence of periodic solutions.  Results that  concern  continuous  and delayed couplings  are considered in \cite{akhmet1, akhmet2}.  

 We believe that the approach proposed in this paper will be useful  for the  investigation of   a wide range of  problems,   focusing not only on  synchrony and  pulse-couplings, but also  phase locking, frequency  locking of systems,  families of oscillators with  continuous couplings.   The method can be used  to analyze  inhibitory models as well as to evaluate the effects of coupling time deviations. Moreover, the model is  suitable for the investigation of the existence of  quasi-periodic and almost  periodic motions.  
\section{The prototype map and two identical oscillators}\label{sec1} 
In this section  we shall define the map, which  is the basic instrument of our investigation. It is constructed for  a  model more general, than  is needed for this paper, to  be the  basis for  future  investigations. 

Let us consider  
two    identical oscillators,  $x_1(t),x_2(t), t \ge 0,$  which  satisfy  the following differential  equations 
 \begin{eqnarray}
 \label{28}
 && x'_i =  f(x_i),
 \end{eqnarray}
 where $ 0 \le x_i \le 1, i = 1,2.$  When the oscillator  $x_j$ fires at  the moment $t$ such  that   $x_j(t) = 1,x_j(t+) = 0,$ then the value of another oscillator with $i \not =j,$  changes   so that  
	 \begin{eqnarray}\label{29}
x_i(t+)=
\left\{ 
\begin{array}{ll} 
0,\:\: {\rm if} \:\:x_i(t) + \epsilon \ge 1, \\ 
x_i(t) + \epsilon,\:\:{\rm otherwise}.
 \end{array}
\right.
\end{eqnarray}
Denote by  $u(t,0,u_0),$ the solution of  the  equation 
\begin{eqnarray}\label{ass}
&& u' =  f(u),
\end{eqnarray}
such  that  $u(0,0,u_0) = u_0.$ 
Assume that   the  solution exists, is unique and continuable to the threshold for all  $u_0.$
Consider the solution $u(t) = u(t,0,v + \epsilon)$ of (\ref{ass}).  Denote by  $s(v)$ the moment when  $u(s) = 1,$ and define  the function $\bar  L(v)  = u(s,0,0)$ on $(0,1-\epsilon).$ 

The following conditions will be needed throughout  the paper:
\begin{itemize}
	\item [(A1)]  $\bar L(v)$  is a  strictly decreasing  continuous  function;
	\item [(A2)]  $\eta = \lim_{v \to 0+} \bar L(v) > 1 - \epsilon;$ 
	 \item [(A3)] $\lim_{v \to  1 -\epsilon} \bar L(v) = 0.$
\end{itemize}
Conditions  $(A1),(A3)$  are valid, if, for example, $f$ is  a positive and lipschitzian function. Another case will be considered in Example \ref{exa4}.  It  is obvious that   there exists a unique fixed point, $v^*, \bar L(v^*)  = v^*.$  

\begin{figure}[hpbt]
  \centering  
  \epsfig {file=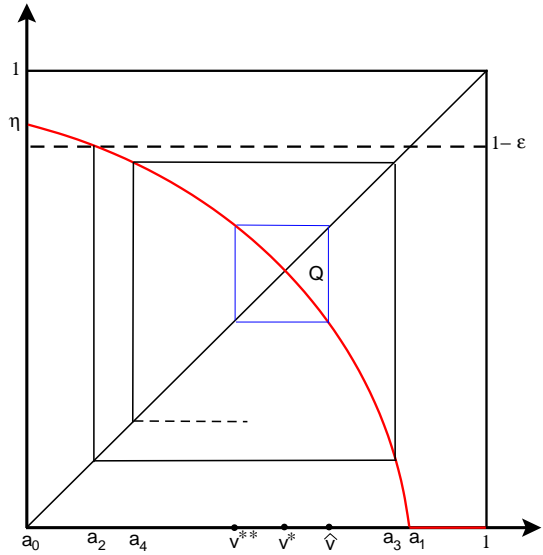, width=2.8in}
  \caption{The graph of  function $w = L(v),$ in red,  and  the   period$-2$ orbit in blue. The  points  $a_0 =0, a_1 =  1-\epsilon, a_{k+1} =  L^{-1}(a_k), k = 1,2,3.$(Color online)} \label{Fig1}
\end{figure}

\begin{figure}[hpbt]
  \centering  
  \epsfig {file=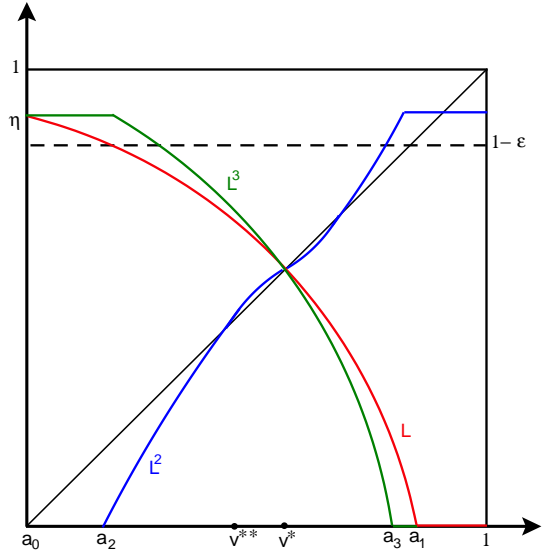, width=2.8in}
  \caption{The graphs of $L,L^2$ and $L^3$  in  red, blue and green respectively.(Color online)} \label{Fig2}
\end{figure}

  Now,  define  a map  $L: [0,1] \to [0,1],$ such that
\begin{eqnarray}\label{31}
L(v) =
\left\{ 
\begin{array}{ll} 
\bar L(v),\:\: {\rm if} \:\: v \in (0, 1-\epsilon), \\
\eta,  \:\: {\rm if} \:\: v = 0,\\
0,\:\:{\rm if} \:\: v \in [1-\epsilon,1].
\end{array}
\right.
\end{eqnarray}	 
 This  newly  defined function  is continuous on $[0,1].$  The sketch  of its   graph is shown  in  Figure \ref{Fig1}.	
	To  make  the following discussion constructive consider  the sequence of maps $L^k(v), k = 1,2,\ldots,$ where $L^k(v) = L(L^{k-1}(v))$ if $k \ge 2.$  Their  graphs with  $k = 1,2,3$ are shown in  Figure \ref{Fig2}.	
	   Denote $a_0 =0, a_1 = 1-\epsilon,  a_2 = L^{-1}(1-\epsilon), a_3 = (L^2)^{-1}(1-\epsilon),\ldots.$  The sequence can be obtained also  through  iterations  $a_0 =0, a_1 = 1-\epsilon, a_{k+1} =  L^{-1}(a_k), k = 1,2,\ldots,$ which are seen in Figure  \ref{Fig1}. 	It  is clear that the sequences $a_{2i}$ and $a_{2i+1}$  are monotonic, increasing and decreasing  respectively.
	    One can  verify   existence of  a fixed point $v^{**} \le v^*$ of the  map $L^2(v)$  such  that  $\hat v= L(v^{**}) \ge v^*,$ and there are no fixed points of  $L^2$ in $(0,v^{**}).$ Moreover,  $a_{2i} \to v^{**} $ and $a_{2i+1} \to \hat v$ as $i \to \infty.$   In the case $v^{**}= v^*,$ there is no non-trivial  period$-2$ points of  $L.$   
	   
	 Let us show  how iterations of $L$  can be  useful for the investigation of  synchronization. 	
	Consider a motion $(x_1(t),x_2(t))$ and  a firing moment  $t_0\ge 0$   such  that $x_1(t_0) = 1, x_1(t_0+) = 0, x_2(t_0+) =v,v \in [0,1].$ 
\begin{lemma} Motion $(x_1(t),x_2(t))$  synchronizes if and only  if  there exists  a number $k$ such that   $1 -\epsilon \le L^k(v) \le 1.$ 
\end{lemma}	 
{\bf Proof}.  Let  us consider only  necessity, since  sufficiency is obvious. We shall consider  the following two  cases: $(\alpha)\,   0 \le v  < 1- \epsilon; \,(\beta)\, 1- \epsilon \le v \le 1.$ 
	 
	  $(\alpha).$ It  is clear that  the couple does not  synchronize at  the moment $t =t_0.$ While it  is not  in  synchrony,  there exists a sequence $t_0 <t_1<\ldots,$ such that $x_1$ fires at  moments  $t_i$ with  even $i$ and $x_2$ fires at $t_i$ with  odd indexes. Set $v_i = x_1(t_i)$ if $i$  is odd, and   
$v_i = x_2(t_i)$ if it  is even.   Use  the definition of $L$ and identity  of oscillators  to  obtain that  $v_{i+1} = L(v_i), i \ge 1.$ This demonstrates that   map $L$  evaluates alternatively  the sequence of values $x_1$ and $x_2$ at firing moments. 
	 
	The pair synchronizes  eventually if and only if there  exists $k \ge 1$ such that $x_1(t) \not = x_2(t),$ if  $t \le t_k,$ and $x_1(t) = x_2(t),$ for $t > t_k.$  Both oscillators  have to fire at $t_k.$ That is,  $1- \epsilon \le v_{k} < 1.$ 
	
$(\beta).$ Consider $1-\epsilon \le v < 1,$ as the case $v=1$ is primitive.   We have that  $t_0$ is a common firing moment of both $x_1$ and $x_2,$ and it  is the synchronization moment. Moreover,   $1-\epsilon < L^2(x_2(t_1)) = \eta < 1.$  The lemma is proved.$\square$

  Thus, it  is   confirmed  that  the analysis of synchronization  is fully  consistent with the dynamics of the introduced map $L(v)$ on $[0,1],$  and,  therefore,  $L$  can  be used  as a valuable tool in further investigation  of the topic.
  
   		Let us consider the rate of synchronization.  We solve the problem by  indicating  initial points which  synchronize after precisely   $k, k \ge 0,$  iterations of the map. Denote by $S_k$ the region in $[0,1],$  where  points $v$ are synchronized after $k$ iterations of  map $L.$ 	One can  see that $S_0 = [1-\epsilon,1], S_1 = [a_0,a_2],$ and  $S_k = (a_{k-1}, a_{k+1}],$ if $k \ge 3,$ is an odd  positive  integer, and   $S_k = [a_{k+1}, a_{k-1}),$ if $k \ge 2,$ is an even  positive integer.
	
	From the discussion made above it follows that the closer  $v$ is  to    $v^{**}$  from the left  or  to  $\hat v$  from the  right, the later is  the moment of synchronization. 	
		
		Denote by $T$ the natural period of  oscillators,  that  is, the period, when  there are no  couplings, and by $\tilde T$  the time needed for  solution $u(t,0,v^*)$  of (\ref{ass}) to  achieve threshold.  Since each oscillator necessarily fires within any  interval of length $T$  and  the distance between two  firing moments of an oscillator  are not   less  than $\tilde T,$ on the basis of the above discussion,  the following assertion is valid.	
	\begin{theorem} \label{thm1} Assume that   conditions $(A1)-(A3)$ are  valid, and   $t_0\ge 0$  is  a firing moment   such  that $x_1(t_0) = 1, x_1(t_0+) = 0.$ If $x_2(t_0+)  \in  S_m, m$  is a  natural  number, then the  couple $x_1,x_2$ synchronizes within  the time interval $[t_0 + \frac{m}{2}\tilde T,t_0 + mT].$
	\end{theorem}	
		  One can  easily  see that whenever condition $(A2)$ is not  valid,
the system  does not  synchronize.
	\begin{example} Consider the   integrate-and-fire model of two  identical  oscillators $x_1,x_2,$  with  the differential equations	
	\begin{eqnarray}
		&& x_i'= x_i^2 +c,	
	\label{25}
	\end{eqnarray}
	where   $i =1,2,$ and $c$ is a positive constant.  It is known  \cite{gerst1} that the  canonical type $I$ phase model \cite{ermen1} can  be reduced by a  transformation  to  the form 
	\begin{eqnarray} \label{be1}
	&& u' = u^2 +c.
	\end{eqnarray}
  This time we investigate  the model with  the pulse-coupling.
  
 Since the two  equations  are  identical, we shall consider a solution  $u(t)$  of  equation (\ref{be1})  to  construct map  $L.$  We have that  $u(t,0, v + \epsilon) = \sqrt{c}\tan(ct + \arctan(\frac{v+\epsilon}{\sqrt{c}}))$ and 
	\begin{eqnarray} \label{e2}
	&& \sqrt{c}\tan(cs + \arctan(\frac{v+\epsilon}{\sqrt{c}})) =1.
	\end{eqnarray}
Next, $u(s,0,0) = \sqrt{c}\tan(cs),$ and by  applying  (\ref{e2}) we find that 
	\[L(v)  = c\frac{1-v-\epsilon}{c+v+\epsilon}, 
\]
if $v \in (0,1-\epsilon),$ and the fixed point  is equal  to  $v^*  = \sqrt{(c +\epsilon /2)^2 + c(1-\epsilon)} - (c+\epsilon /2).$  	
	Evaluate	
	\[L(0)  = c\frac{1-\epsilon}{c+\epsilon}
\]
 to  see that  $L(0) < 1 - \epsilon,$ and condition $(A2)$ is not  valid.
Moreover,  one can  verify  that  $L'(v) < 0.$ 
  
    Thus,  we obtain that   the couple  does not  synchronize, and our simulations
confirm this.  
	\end{example}	
	\begin{example} \label{pesk}  Consider the  following integrate-and-fire model  of  two identical  oscillators,
	  $x_1,x_2,$  with  the differential equations	
	\begin{eqnarray}
		&& x_i'= S -  \gamma x_i,	
	\label{e4}
	\end{eqnarray}
	where $i=1,2,$ positive constants  $S,\gamma$ satisfy $\kappa = \frac{S}{\gamma} >1.$ 
	   One can  find that $u(t,0, v + \epsilon) =
	(v  +\epsilon) {\rm e}^{-\gamma t} + \kappa (1-{\rm e}^{-\gamma t})$ and 	$u(s,0,0) = \kappa (1-{\rm e}^{-\gamma s}).$ The last  two expressions  imply  that  
\begin{eqnarray}
	&& L(v) = \kappa \frac{1- (v +\epsilon)}{\kappa -  (v +\epsilon)},	
	\label{e6}
\end{eqnarray}
if $0 < v < 1-\epsilon.$ 
 
 There is  a unique  fixed point of $L$ and $L^2,$ and it is equal to 
  \begin{eqnarray}
	&&  v^*  = (\kappa  - \frac{\epsilon}{2}) -  \sqrt{\kappa^2 - \kappa + \frac{\epsilon^2}{4}}. 
	\label{13}
	\end{eqnarray}
	Finally,  $L(0) = \kappa \frac{1- \epsilon}{\kappa -  \epsilon} > 1- \epsilon.$  That  is, all conditions of the last  theorem are valid, and   the assertion in \cite{p} is  proved. 
	\end{example}
		\begin{remark} Map  $L$  is similar to  that  in \cite{p}, but the  argument here is a  coordinate before a jump, while in  the  paper  the  argument is a coordinate  after a  jump.  This difference  is not   a critical one.    The most important  point   is that  C. Peskin uses it  only  as an auxiliary device to build the Poincar\'e map.  We use   $L$  itself  as the main map, with a newly  defined continuous extension, which  simplifies  the discussion in this section and throughout our investigation.  We  revisit  the problem of two identical oscillators, since $L$ is the prototype map  in  our  analysis. In addition to  the main synchronization result,  regions with equal time of synchronization are  indicated, and the  value of the  fixed point, $v^*,$  is evaluated. 
\end{remark} 
		\begin{example}\label{exa4} Consider the following   integrate-and-fire system of pulse-coupled oscillators, $x_1,x_2,$ such  that
		\begin{eqnarray}
	&& x_1'= f(x_1), \nonumber\\  
	&& x_2'= f(x_2),	
	\label{152}
	\end{eqnarray}  
	where 	
	\begin{eqnarray}
f(s) =
\left\{ 
\begin{array}{ll} 
4 - 3s , & {\rm if}\:\: 0 < s \le  1/3; \\ 
3, & {\rm if}\:\: 1/3 < s \le 2/3;\\
4 - 3(s-2/3) , & {\rm if}\:\: 2/3  < s \le  1. 
\end{array}
\right.
\end{eqnarray}	
We have found that   map  $L$ for  this system exists and is equal to 
	\begin{eqnarray}
L(v) =
\left\{ 
\begin{array}{ll} 
2\frac{2-3v-3\epsilon} {4-3v-3\epsilon} , & {\rm if}\:\: 0 < v \le  1/3 - \epsilon; \\ 
1-v-\epsilon, & {\rm if}\:\: 1/3-\epsilon < v \le 2/3 -\epsilon;\\
\frac{4}{3}\frac{1-v-\epsilon} {2-v-\epsilon}, & {\rm if}\:\: 2/3-\epsilon  < v \le  1. 
\end{array}
\right.
\end{eqnarray}	
One can check that   conditions $(A1)-(A3)$ are fulfilled for this map. Moreover,   the fixed points are  equal  to  $v^* = (1-\epsilon)/2, v^{**} = 1/3$ and $\hat v = 2/3 -\epsilon.$ Finally,  all the motions, which  start outside of the periodic trajectory  synchronize eventually, and all of them  are periodic inside the trajectory. 
	\end{example} 	
The last  example shows that  the assumptions for   map $L,$ including the  existence of the non-trivial period$-2$ motion,  can  be realized for  even the differential equations with  discontinuous right-hand side.  Moreover,  in  future investigations one can  consider  isolated periodic solutions, stable or unstable.   The theoretical  consequences of this research are clear, if one uses  the mappings theory, but construction of examples requires    additional time. 
	\begin{example} \label{exa1} Consider the  model  of two integrate-and-fire identical  oscillators, $x_1,x_2,$  which are  pulse-coupled and  
	 \begin{eqnarray}
	&& x_1'= S - \gamma x_1 + \beta x_2, \nonumber\\  
	&& x_2'= S - \gamma x_2 + \beta x_1,	
	\label{15}
	\end{eqnarray}
where constants $S,\gamma$ and $\beta$ are positive numbers. One can  easily  see that  the system is the extended  Peskin's model in Example \ref{pesk}. The terms with  coefficient $\beta$ are newly introduced in the system. They  reflect  the permanent  influence  of the partners  during the process. Eigenvalues  associated   to  (\ref{15}) are  $\lambda_1 = 	-\gamma + \beta$ and  $\lambda_2 = 	-\gamma - \beta.$ We suppose that  $\beta$ is  small  so  that  both eigenvalues are negative. Moreover,  it is assumed  $\kappa = S/\gamma > 1.$  Then,
 $\kappa_1 = -S/\lambda_1 >1$  if $\beta$ is sufficiently small.  
  The solution of system (\ref{15}) with  value $(0,v+\epsilon)$ at  $t = 0,$ is equal to 
	$u_1(t)= \frac{1}{2}[{\rm e}^{\lambda_1 t} - {\rm e}^{\lambda_2 t}](v+\epsilon) - \kappa_1({\rm e}^{\lambda_1 t}-1),\,u_2(t)= \frac{1}{2}[{\rm e}^{\lambda_1 t} + {\rm e}^{\lambda_2 t}](v+\epsilon)  - \kappa_1({\rm e}^{\lambda_1 t}-1).$
By  using these expressions   obtain  the equation $ \frac{1}{2}[{\rm e}^{\lambda_1 s} + {\rm e}^{\lambda_2 s}](v+\epsilon) + \kappa_1(1-{\rm e}^{\lambda_1 s}) = 1,$
 and 	  construct  
	$L(v) = 1 - (v+\epsilon){\rm e}^{\lambda_2 s}.$
	Map $L$ is too  complex  to  analyze  for properties $(A1)-(A3).$ That  is why  we  will compare this model with   the couple in Example \ref{pesk}. The last  two  equations imply  $L(0) = \kappa\frac{1  -\epsilon}{\kappa -\epsilon}  > 1-\epsilon,$  if  $\beta =0$ and $v= 0.$    That  is, if $\beta$ is sufficiently small, then   condition $(A2)$ is valid. We have found, also, by  direct evaluations that the derivative $L'(v)$
is  negative if $S$ and $\beta$ are sufficiently  large and small respectively. That  is,   condition $(A1)$  is  fulfilled. It  is easy to  verify that  condition $(A3)$ is also  correct.	
	
	 Now, using the continuity  theorem in parameters \cite{hartman},  one can  find  that  map $L$  may  admit a period$-2$ point only  if the orbit is  as close to  the fixed point $v^*$ of  Example \ref{pesk}  as  $\beta$  is  small. Consequently,  the measure of the set  of  points, which  can  not  be synchronized diminishes  as $\beta \to 0.$   This result is a new one.  In previous papers the differential equations were separated.  		
		\end{example}	
\section{The multidimensional system of non-identical oscillators.}  
   Consider the  model of $n$ non-identical oscillators   given by  relations (\ref{1})  and (\ref{2}). The domain  of this model consists of points $x  = (x_1,x_2,\ldots,x_n)$ such  that $0 \le x_i  \le 1 +\zeta_i(x)$ for all $i =1,2,\ldots,n.$ 
    
   Fix two of the  considered oscillators, let  say,  $x_l$ and $x_r.$     
\begin{lemma}\label{lem} Assume that  condition $(A2)$ is  valid, and   $t_0\ge 0$  is  a firing moment   such  that   $x_l(t_0) = 1  + \zeta_l(x(t_0)), x_l(t_0+) = 0.$ If   parameters are sufficiently close to  zero,  and  absolute values  of   parameters of perturbation   sufficiently small with  respect to  $\epsilon,$  then the couple $x_l, x_r$ synchronizes within the time interval $[t_0, t_0 + T]$ if  $x_r(t_0+) \not \in [a_0,a_{1})$  and  within the time interval  $[t_0 + \frac{m-1}{2}\tilde T, t_0 + (m+1)T],$ if  $x_r(t_0+) \in S_{m},  m \ge 1.$ 
\end{lemma}
{\bf Proof}.   Denote by $x(t) = (x_1(t),x_2(t),\ldots, x_n(t)),$ the  motion of the oscillator. 
   If  $1 + \zeta_r(x(t_0)) -\varepsilon - \varepsilon_{r} \le x_r(t_0)\le  1 +\zeta_r(x(t_0)),$ then these two oscillators fire simultaneously, and we   only need  to  prove the  persistence of  synchrony  which   will be done later. So, fix another  oscillator $x_r(t)$ such  that $ 0 \le  x_r(t_0) < 1 +\zeta_r(x(t_0)) -\varepsilon - \varepsilon_{r}.$   
   
  While the pair does not synchronize, there  exists a sequence of moments $0<t_0<t_1<\ldots,$ such  that  oscillator $x_l$ fires  at $t_i$  with  even  $i,$  and $x_r$ fires at $t_i$ with  odd  $i.$   For the sake of brevity  let $u_i = x_l(t_i), i = 2j +1,  u_i = x_r(t_i), i = 2j, j \ge 0.$  In what follows we shall  evaluate  the difference $u_{i+1} - L(u_{i}).$ 

 Let us fix an even $i$ and $u_i  = x_r(t_i).$  If the parameters are sufficiently  small, then there are $k \le n-2$ distinct firing moments of the motion   $x(t)$  on  the interval $(t_i,t_{i+1}).$ Denote by   $t_i <\theta_1<\theta_2<\ldots<\theta_k <t_{i+1},$   the moments  of firing, when at least one of the coordinates of $x(t)$ fires, and   $v(t,t_0,v_0)$ the solution of the  equation  (\ref{1})  with  $v(t_0,t_0,v_0)=v_0.$    We have that 
  \begin{eqnarray} \label{4g3}
  &&x_r(\theta_1)  =  x_r(t_i) + \epsilon + \int_{t_i}^{\theta_1}f(x_r(s))ds +\nonumber\\
  &&\int_{t_i}^{\theta_1}\phi_r(x(s))ds,
  \end{eqnarray}
 where $x(t) = v(t,t_i,x(t_i+)),$  
   
  \begin{eqnarray} \label{4g31}
  &&x_r(\theta_2)  =  x_r(\theta_1) + \epsilon + \int_{\theta_1}^{\theta_2}f(x_r(s))ds + \nonumber\\
  &&\int_{\theta_1}^{\theta_2}\phi_r(x(s))ds,
   \end{eqnarray}
   where $x(t) = v(t,\theta_1,x(\theta_1+)),$  
     
	\[\ldots\ldots\ldots\ldots\ldots\ldots\ldots\ldots\ldots\]
  \begin{eqnarray} \label{4g32} 
 &&x_r(t_{i+1}) =   x_r(\theta_k) + \epsilon + \int_{\theta_k}^{t_{i+1}}f(x_r(s))ds +\nonumber\\
 &&\int_{\theta_k}^{t_{i+1}}\phi_r(x(s))ds,
\end{eqnarray} 
where $x(t) = v(t,\theta_1,x(\theta_k+)).$ 

The moment $t_{i+1}$  satisfies 
\begin{eqnarray} \label{4g5}
1 +\zeta_r(x(t_{i+1})) - \epsilon -\epsilon_r  \le x_r(t_{i+1}) \le 1 +\zeta_r(x(t_{i+1})).
\end{eqnarray} 
Similarly  to  the expressions for $x_r$  one can  obtain  
\begin{eqnarray} \label{4g4}
&&x_l(\theta_1)  = \int_{t_i}^{\theta_1}f(x_l(s))ds +\int_{t_i}^{\theta_1}\phi_l(x(s))ds, \nonumber\\
  &&x_l(\theta_2)  =   x_l(\theta_1) + \epsilon + \int_{\theta_1}^{\theta_2}f(x_l(s))ds +\nonumber\\
  &&\int_{\theta_1}^{\theta_2}\phi_l(x(s))ds, \nonumber\\
  && \nonumber\\
  &&\ldots\ldots\ldots \ldots\ldots\ldots \ldots\ldots\ldots\nonumber\\
   &&   \nonumber\\
 &&x_l(t_{i+1}) =   x_l(\theta_k) + \epsilon + \int_{\theta_k}^{t_{i+1}}f(x_l(s))ds + \nonumber\\
 &&\int_{\theta_k}^{t_{i+1}}\phi_l(x(s))ds.
\end{eqnarray} 

Formulas  (\ref{4g3})  to   (\ref{4g4})   define  $u_{i+1}= x_l(t_{i+1}).$ Similarly one can  evaluate the number for  odd $i.$ 

Let  us now  find the value of $L(u_{i}).$ 
With  this aim,  evaluate   
 \begin{eqnarray} \label{5g3}
&& \phi(\bar t_{i+1}) = x_r(t_i) + \epsilon + \int_{t_i}^{\bar t_{i+1}}f(\phi(s))ds ,
\end{eqnarray} 
where $\bar t_{i+1}$ satisfies $\phi(\bar t_{i+1}) =1,$
 and 
\begin{eqnarray} \label{5g4}
&&\psi(\bar t_{i+1})  = \int_{t_i}^{\bar t_{i+1}}f(\psi(s))ds ,
\end{eqnarray} 
to   find that    $L(u_{i}) = \psi(\bar t_{i+1}).$  Next, we  will show that   the difference  $u_{i+1} - L(u_{i})$ is small if the parameters are small.  

First,  one can  find   
 \begin{eqnarray} \label{4g33}
  && \phi(t) - x_r(t)  =  \nonumber\\
  &&\int_{t_i}^{t}[f(\phi(s))- f(x_r(s))]ds - \int_{t_i}^{t}\phi_r(x(s))ds,
   \end{eqnarray}
   for $t  \in [t_i,\theta_1].$ 
   
   Then, by  applying the Gronwall-Bellman Lemma one can   easily  see  
   \begin{eqnarray} \label{4g34}
  && |\phi(\theta_1) - x_r(\theta_1)|  \le \mu_r (\theta_1 - t_i) {\rm e}^{\ell(\theta_1 - t_i)}, 
   \end{eqnarray}
    where  $\ell$ is the Lipschitz constant of $f.$ 
    Next, we have 
     \begin{eqnarray} \label{4g35}
  && |\phi(\theta_2) - x_r(\theta_2)| \le  \nonumber\\
  &&[\mu_r(\theta_1 - t_i){\rm e}^{\ell(\theta_1 - t_i)} +\mu_r(\theta_2 - \theta_1)+ \epsilon]{\rm e}^{\ell(\theta_2-\theta_1)}, 
   \end{eqnarray}  
    if $t \in [\theta_1,\theta_2].$
    
Without loss of generality, assume that  $t_{i+1} > \bar t_{i+1}.$  Proceeding the  evaluations made above,  we can  obtain   $|1 - x_r(\bar t_{i+1})| =  |\phi(\bar t_{i+1}) - x_r(\bar t_{i+1})| =  \Phi_1(\epsilon,\mu_r),$   where 
	\[\Phi_1(\epsilon,\mu_r) \equiv \mu_r [(\theta_1 - t_i){\rm e}^{\ell(\bar t_{i+1} - t_i)} + \sum_{j=1}^{k-1}(\theta_{j+1} -\theta_{j}) {\rm e}^{\ell(\bar t_{i+1} - \theta_j)} + \]\[(\bar t_{i+1} - \theta_k){\rm e}^{\ell(\bar t_{i+1} - \theta_k)}] +  \epsilon \sum_{j=1}^{k} {\rm e}^{\ell(\bar t_{i+1} - \theta_j)}.\] There are  positive  numbers $\mu$  and $M,$ which  satisfy  $\mu \le f(s) \le M,$ if $0 \le s \le 1 + \max_{i} \xi_i.$ 
One  can  request  the following inequality  $\max_{i= 1,\ldots,n} \mu_i < \mu.$	
  We have that  $ |x_r(t_{i+1}) -  x_r(\bar t_{i+1})| \le |1 - x_r(t_{i+1})| + |1 - x_r(\bar t_{i+1})| \le  \Phi_1(\epsilon,\mu_r) + \xi_r.$ Consequently,  
	\[|t_{i+1} - \bar t_{i+1}| < \frac{\Phi_1(\epsilon,\mu_r) + \xi_r}{\mu - \mu_r} \equiv \Phi_2(\epsilon,\mu_r,\xi_r).
\]
By  applying   (\ref{4g4})  and (\ref{5g4}),  making similar  evaluations for  (\ref{4g34}) and (\ref{4g35}) one can  find 
$|\psi(\bar t_{i+1}) - x_l(\bar t_{i+1})| \le \Phi_1(\epsilon,\mu_l).$

 Then,  we  have that $|u_{i+1} - L(u_{i})| = |\psi(\bar t_{i+1}) - x_l( t_{i+1})| \le  |\psi(\bar t_{i+1}) - x_l( \bar t_{i+1})|+ |x_l( t_{i+1}) - x_l( \bar t_{i+1})|,$ and, consequently,
 \begin{eqnarray}\label{last}
 &&  |u_{i+1} - L(u_{i})|\le  \Phi(\epsilon,\mu_r,\mu_l,\xi_r),
   \end{eqnarray} 
 where $\Phi \equiv \Phi_1 + \Phi_2(M + \mu_r).$  
 It  is obvious that  $\Phi$   tends to  zero  as  the parameters  do.  
This convergence is uniform  with  respect to  $u_0.$ We can also  vary the number of points $\theta_i$ and their location  in the intervals $(t_j,t_{j+1})$ between $0$ and $n-1.$ The  convergence  also is indifferent with  respect to  these variations. 
  
Consider the sequence of inequalities
$$|u_i - L^i(u_0)| \le |u_i - L(u_{i-1})| + | L(u_{i-1}) -L(L^{i-1}(u_0))|, i = 1,2,\ldots.$$ Then recurrently, by  applying continuity  of  $L,$  (\ref{last})  and   $L^m(u_0) \in [1-\epsilon,1],$   conclude that either $1+\xi_r-\epsilon-\epsilon_r \le u_m < 1+\xi_r$ or  $1+\xi_l-\epsilon-\epsilon_l \le u_{m+1} < 1+\xi_l,$ if the parameters are sufficiently  close to zero,  and  absolute values of the parameters of perturbation  are sufficiently small with  respect to  $\epsilon.$  Both of these inequalities bring the pair  to  synchronization. 
 
 Since each of the iterations of map $L$ happens within an interval with length not  more than $T,$ we obtain that  couple $x_l,x_r$  synchronizes no later than  $t= t_0 + (m+1)T.$   Similarly,  the couple synchronizes not  earlier than 
 $ t = t_0 +  \frac{m-1}{2}\tilde T.$
   
 If two oscillators $x_l$ and $x_r$ are non-identical  and fire simultaneously at  a moment $t = \theta,$ how  will  they retain the state of firing in unison, despite being different?  To find  the required conditions, let  us  denote by   $\tau, \tau > \theta$ a moment  when  one of them, let's  say $x_r,$ fires. We have that $x_l(\theta+)  =x_r(\theta+) =0.$  Then  $x_l(t)  =x_r(t), \theta \le t \le \tau.$  It is clear that to  satisfy $x_l(\tau+)  =x_r(\tau+) =0,$ we need 
$1 + \zeta_l(x(\tau)) - \epsilon - \epsilon_l \le x_l(\tau).$ By applying formula  (\ref{4g5}) again, this time with  $t_i = \theta, t_{i+1} = \tau,$ one can  easily  obtain that  the inequality  is correct if  parameters   are close to  zero,  and  absolute values of the parameters of perturbation  are sufficiently small with  respect to  $\epsilon.$  Thus,  one can  conclude that  if  a couple of oscillators is synchronized at   some moment of time then it persistently continues to  fire in unison. The lemma is proved.$\square$
\begin{remark} The last  lemma not only plays an auxiliary  role for next main theorem, but  can also be considered  a  synchronization result for  the model of two non-identical oscillators.
\end{remark}
 Let us  extend the result of the lemma  for the whole ensemble.  	 
\begin{theorem}\label{thm}  Assume that   condition $(A2)$ is  valid,  and   $t_0\ge 0$  is  a firing moment   such  that   $x_j(t_0) = 1  + \zeta_j(x(t_0)), x_j(t_0+) = 0.$  If the  parameters  are sufficiently close to  zero,  and  absolute values of  parameters of perturbation  are sufficiently small with  respect to  $\epsilon,$  then   the motion   $x(t)$ of the system  synchronizes  within the time  interval  $[t_0,t_0 +T],$ 
if   $ x_i(t_0+) \not \in [a_0,a_1), i \ne j,$   and  within the time interval  $[t_0 + \frac{\max_{i \ne j}k_i-1}{2}\tilde T,t_0 +(\max_{i \ne j}k_i+1)T],$  if  there exist  $x_s(t_0+) \in [a_0,a_1)$ for some $s \not = j$  and  $x_i(t_0+) \in S_{k_i}, i \not = j.$
\end{theorem} 
{\bf Proof.}  Consider the non-trivial case. Applying the  last  lemma we can see that each pair $(x_j,x_i),  i \ne j,$ synchronizes  within  $[t_0 + \frac{\max_{i \ne j}k_i-1}{2}\tilde T, t_0 + (\max_{i \ne j}k_i+1)T].$  The  theorem is proved. $\square$

Now,  replace  coupling (\ref{2})  by  
	 \begin{eqnarray}\label{2'}
x_i(t+)=
\left\{ 
\begin{array}{ll} 
0,\:\: {\rm if} \:\:x_i(t) + \epsilon +\epsilon_{i} \ge 1+\zeta_i(x), \\ 
x_i(t) + \bar \epsilon +\epsilon_{i},\:\:{\rm otherwise}, \end{array}
\right.
\end{eqnarray}
where $\bar \epsilon, \bar \epsilon + \epsilon_{i} > 0,$ is a new parameter, independent of $\epsilon.$ Consider system (\ref{1})  with   (\ref{2'}). One can  find that  the following assertion is valid.  
\begin{theorem}\label{thmnew}  Assume that   condition $(A2)$ is  valid,  and   $t_0\ge 0$  is  a firing moment   such  that   $x_j(t_0) = 1  + \zeta_j(x(t_0)), x_j(t_0+) = 0.$  If   parameters  $\bar \epsilon,\mu_i, \xi_i, \epsilon_{i},$ are sufficiently close to  zero,  then   the motion   $x(t)$ of the system  synchronizes  within the time  interval  $[t_0,t_0 +T],$ 
if   $ x_i(t_0+) \not \in [a_0,a_1), i \ne j,$   and  within the time interval  $[t_0 + \frac{\max_{i \ne j}k_i-1}{2}\tilde T,t_0 +(\max_{i \ne j}k_i+1)T],$  if  there exist  $x_s(t_0+) \in [a_0,a_1)$ for some $s \not = j$  and  $x_i(t_0+) \in S_{k_i}, i \not = j.$
 \end{theorem}
 
We can  see that   (\ref{2'}) changes the style of interaction  in the model. It  depends  on  distance of oscillators to  thresholds.  We use this  to  introduce  delay and continuous couplings in papers \cite{akhmet1}  and  \cite{akhmet2}, respectively.

\begin{figure}[hpbt]
  \centering  
  \epsfig {file=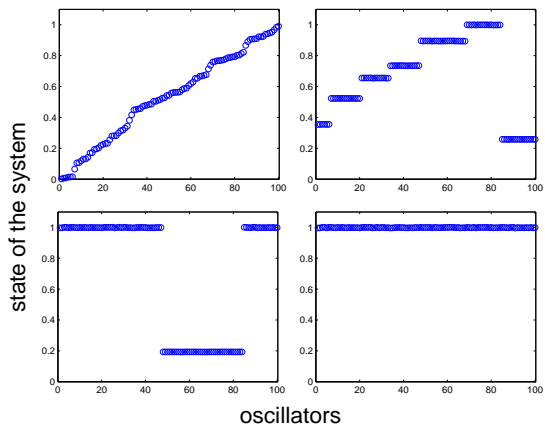, width=3.4in}
  \caption{The  state of the model   before the first, twenty  first, forty  second and sixty third  firing of the system.  The flat  sections of the graph are groups of  synchronized oscillators.(Color online)} \label{Fig6}
\end{figure}

To  illustrate Theorem \ref{thm},  consider a  group  of  oscillators, $x_i, i = 1,2,\ldots,100,$  with random  uniform distributed  start values in  $[0,1].$ It  is supposed that they satisfy  the  equations $x_i' =  (3 +0.01\bar\mu_i)  - (2 +0.01\bar\zeta_i)x_i.$ The constants  $\bar\mu_i, \bar\zeta_i,$  as well as $\bar\xi_i$ in the  thresholds  $1 +  0.005\bar\xi_i, i = 1,2,\ldots,100,$ are uniform random distributed  numbers from  $[0,1].$  In Figure \ref{Fig6} one can see the result of simulation with  $\epsilon = 0.08,$ where the  state of the system is shown  before  the first,   twenty  first, forty  second and sixty third  firing of the system. So, it is obvious that  eventually the model  shows  synchrony. 
Let us  describe  a more general system of  oscillators  such  that  Theorem \ref{thm} is still true. 
  A system of $n$ oscillators is given, such that  if $i-$th oscillator does not fire or jump up, then it satisfies the $i-$th equation  of system (\ref{1}). If  several oscillators $x_{i_s}, s = 1,2,\ldots,k,$ fire  so  that $x_{i_s}(t) = 1 +\zeta_{i_s}(x)$ and  $x_{i_s}(t+) = 0,$ then all other oscillators $x_{i_p}, p= k+1,k+1,\ldots,n,$ change their coordinates by the law
 \[x_{i_p}(t+)=
\left\{ 
\begin{array}{ll} 
0,\:\: {\rm if} \:\:x_{i_p}(t) + \epsilon + \sum_{s=1}^k \epsilon_{i_pi_s} \ge 1 + \zeta_{i_p}(x), \\ 
x_{i_p}(t) + \epsilon + \sum_{s=1}^k \epsilon_{i_pi_s},\:\:{\rm otherwise}.  \end{array}
\right.\]
One can  easily  see that the last  theorem is correct  for the model  just  described, if 
$\epsilon + \sum_{s=1}^k \epsilon_{i_pi_s} > 0,$ for all possible $k,i_p$ and $i_s.$  
\begin{remark} The analysis  of non-identical oscillators  with non-small  parameters  may  shed light on the investigation of arrhythmias,   chaotic flashing of fireflies, etc. Namely,  the dynamics  in  the neighborhood  and inside the periodic trajectory, $Q$ in  Figure \ref{Fig1},   can be very  complex. We do not exclude the possibility  of  chaos  and  fractals \cite{timme1}. Bifurcation of periodic solutions  can be discussed, if the parameters are small.
\end{remark}
 \section{The Kamke condition and synchronization} \label{Kamke}
In this section we  consider  an  integrate-and-fire  model  with a new type of  continuous connection.  

 We believe that models  of identical  oscillators  with  more general
differential  equations, 
\begin{eqnarray}
 \label{32}
 && x'_1 =  g(x_1,x_2,\ldots,x_n),\nonumber\\
 && x'_2 =  g(x_2,x_3,\ldots,x_1),\nonumber\\
 &&\ldots \ldots \ldots\nonumber\\
 && x'_n =  g(x_n,x_1,\ldots,x_{n-1}), 
 \end{eqnarray}
where $0\le x_i \le 1, i = 1,2,\ldots,n,$  are of both theoretical and applied
interest.   The  positive valued function $g(y_1,y_2,\ldots,y_n)$ in (\ref{32}) is continuously  differentiable  and  indifferent  with  respect  to  permutations of coordinates $y_2$ to  $y_n.$  
 
  When the oscillator  $x_j$ fires at  the moment $t$ such  that   $x_j(t) = 1,x_j(t+) = 0,$ then the value of an  oscillator with $i \not =j,$  changes   so that  
	 \begin{eqnarray}\label{33}
x_i(t+)=
\left\{ 
\begin{array}{ll} 
0,\:\: {\rm if} \:\:x_i(t) + \epsilon \ge 1, \\ 
x_i(t) + \epsilon,\:\:{\rm otherwise}.
 \end{array}
\right.
\end{eqnarray}
  Consider the  cone $\mathbb R^n_+ \subset \mathbb R^n$ of all vectors with  nonnegative coordinates.  Introduce a partial  order  in the cone such that  $a \le b$ if $a_i \le b_i, i = 1,2,\ldots,n,$ \cite{smith}. 
  
  We say  that  function $g$ is of type $\cal K$ in $\mathbb R^n_+$ if $a \le b$ implies that  $g(a) \le g(b).$ The sufficient condition for that  is $\frac{\partial g(y)}{\partial y_i} \ge 0, i \ne 1.$ 
  
  Let  $u(t,t_0,u_0)$ and  $u(t,t_0,u_1), u_0,u_1 \in \mathbb R^n_+,$ be  solutions of (\ref{32}).
  If $g$ is of type $\cal K,$ then \cite{kamke,smith} the dynamics of (\ref{32}) is monotone for $t \ge 0.$ That  is,  $u(t,0,u_0) \le u(t,0,u_1),$ if $u_0 \le u_1.$

 Consider first  the model  of two  oscillators.  Define  map  $L$ for this system in the following way. Take 
the solution $u(t) = u(t,0,(0,v + \epsilon)) =  (u_1,u_2).$ Denote by  $s(v)$ the moment when  $u_2(s) = 1,$ and define  the function $\bar  L(v)  = u_1(s)$ on $(0,1-\epsilon).$  Then define  map  $L$  through (\ref{31}).  Let  us check if conditions $(A1),(A3)$ are valid for this map. Indeed, the continuity  of $\bar L$ is obvious. It  is non-increasing since the monotonicity. Assume that  there exist  numbers $v_1,v_2 \in  (0,1-\epsilon)$ such  that  $v_1 < v_2$ and  $s =  s(v_1) = s(v_2).$ Then, we have a contradiction as the open interval $(v_1,v_2)$ is mapped to  the closed set $\{1\}.$  That  is, $\bar L$ satisfies $(A1).$

Condition   $(A3)$  is easily  verifiable.    Now, one can  determine  sets $S_i,$ similarly  to  that  in Section \ref{sec1}, and prove that  the following theorem  is valid. 
	\begin{theorem} \label{thm4} Assume that $g$ is of $\cal K$ type,   $(A2)$ is  valid,  and   $t_0\ge 0$  is  a firing moment   such  that   $x_1(t_0) = 1, x_1(t_0+) = 0.$ 
If $x_2(t_0+)  \in  S_m, m$  is a  natural  number, then the  couple $x_1,x_2$ synchronizes within  the time interval $[t_0 + \frac{m}{2}\tilde T,t_0 + mT].$
	\end{theorem}
  We have, moreover, that, if $g$ is of $\cal K$ type and  $(A2)$ is not  true then the system does not synchronize.
    
 Consider the multidimensional system of oscillators.  Introduce the function $G(y,z)  \equiv g(y,z,z,\ldots,z),$ and define the integrate-and-fire model of two identical  oscillators $y$ and $z$  with the following system of differential equations
  \begin{eqnarray}
 \label{34}
 && y' = G(y,z),\nonumber\\
 && z' = G(y,z).
 \end{eqnarray} 
  Denote by  $u = (y,z), u(t) = u(t,0,(0,v + \epsilon)),$ the solution of (\ref{34}), and  by  $s(v)$ the moment when  $z(s) = 1.$ Next,  define  the function $\bar  L(v)  = y(s)$ on $(0,1-\epsilon).$  Then  map  $L$  can  be defined by  (\ref{31})  as well as correspond sets $S_i.$  
  By  applying the monotonicity  of the dynamics, one can  prove the following assertion, in a way very  similar to that of   Theorem \ref{thm}.   
\begin{theorem}\label{thm5}  Assume that  $0 \le  \frac{\partial g(y)}{\partial y_i} < \eta, i  \ne 1,$  condition $(A2)$ is  valid, and
   $t_0\ge 0$  is  a firing moment   such  that   $x_j(t_0) = 1, x_j(t_0+) = 0.$     If parameter  $\eta$ is sufficiently  small  then   the motion   $x(t)$ of the system  synchronizes  within the time  interval  $[t_0,t_0 +T],$ if   $ x_i(t_0+)  \in  S_0, i \ne j,$   and  within the time interval  $[t_0 + \frac{\max_{i \ne j}k_i-1}{2}\tilde T,t_0 +(\max_{i \ne j}k_i+1)T],$  if  there exist  $x_s(t_0+) \in [a_0,a_1)$ for some $s \not = j$  and  $x_i(t_0+) \in S_{k_i}, i \not = j.$
\end{theorem}   
  Example \ref{exa1} shows  that the analysis of  the map  $L$ for Theorems \ref{thm4} and \ref{thm5} is not  simple   even with  linear  differential equations. The  results of this section  are therefore are provided  for  numerical application, as well as  for   future investigations. 
 \section{Conclusion} A version of the integrate-and-fire model of pulse-coupled  and non-identical oscillators is investigated in  this  paper.   We  have made significant advances for the solution  of  second
Peskin's conjecture, though we have not  showed that  the measure of
non-synchronized initial  values is zero  as  it was shown in \cite{ms}
for identical oscillators. However, we have  located  non-synchronized  points, and  showed that  the time of synchronization infinitely increases  as the region  of points  to be  synchronized  enlarges. Moreover, it is not  necessary,  in applications, for all points of the domain  to fire in unison, and it is sufficient for the neighborhood of a motion to be synchronized. 

 A prototype map  is introduced that  helps to  precise the results for a system of two identical oscillators and to  solve the problem for the  multi-non-identical-oscillators model.  The approach  of the paper is  universal. For example,  it  can  be applied  if the thresholds  are of the form  $1+ \phi(t,x)$ or  $1+ \phi(t),$ which  represents an oscillating signal  in physiology \cite{glass2} and a variation of the threshold  of the electronic relaxation oscillator \cite{pik}.  Moreover, the jump's  value  $\epsilon$ as well as its perturbations may depend on  $x.$  The method  can  be easily  extended  for models, where   differential  equations have discontinuous right-hand sides,  and only the  existence of the  map  $L$  with  request  properties are important  to  get  appropriate results. Not  only  exhibitory, but  inhibitory  models of oscillators can  be investigated, as well as  delay  of  couplings \cite{buck1,epg,gerst,akhmet1}, and continuous couplings generated by   firing \cite{van,akhmet2}.  One can  consider the models of our paper as discontinuous  {\it cooperative} systems \cite{smith}. Consequently, we expect that  by  applying  methods of  dynamical systems  with  variable moments of discontinuity  \cite{akhmet}, the results  on monotone systems \cite{smith, hirsh,hirsh1} can be extended for  these models. 

 There  is a  rich  collection of results on   oscillators,  obtained through experiments and simulations. The approach  of the present paper can  give theoretical background for them and also  form a basis  for  new  ones. It  can  be applied   not only to the  problems of synchronization, but also to periodic, almost  periodic motions, and the complex behavior of biological models. New small-world phenomena can be discovered.  One can  now request  a certain property  for  the map  $L$ and then  look  for a system which  meets the property. Thus, many new  theoretical challenges can  be brought under  discussion.  Conversely, if  a  system is given, then one can  construct  the corresponding map $L,$ and by  analyzing    find  new features which  have not been  mentioned in the present investigation.  
 \begin{Ack}The author thanks  Mehmet Turan and Zahid Samancioglu  for   the technical assistance.
\end{Ack}
 
\section{References}

\bibliographystyle{model1-num-names}
\bibliography{elsarticle-num.bst}

\end{document}